\documentclass[11pt]{amsart}

\usepackage{amsfonts}

\usepackage{mathrsfs}
\usepackage{bbm}
\usepackage{amsmath}
\usepackage{amssymb}
\usepackage{amscd}
\usepackage{amsthm}

\usepackage{amsaddr}

\makeatletter
\@namedef{subjclassname@2020}{\textup{2020} Mathematics Subject Classification}
\makeatother

\newtheorem{thm}{Theorem}[section]

\newtheorem{cor}[thm]{Corollary}

\numberwithin{equation}{section}

\newcommand{\di}{\displaystyle}
\newcommand{\te}{\textstyle}

\newcommand{\bc}{\mathbb{C}}
\newcommand{\bff}{\mathbb{F}}

\newcommand{\bn}{\mathbb{N}}

\newcommand{\br}{\mathbb{R}}
\newcommand{\bz}{\mathbb{Z}}
\newcommand{\bt}{\mathbb{T}}

\newcommand{\fg}{\mathfrak{g}}

\newcommand{\cf}{\mathcal{F}}

\newcommand{\sce}{\mathscr{E}}

\newcommand{\autg}{Aut\left(G \right)}
\newcommand{\autag}{Aut_a\left(G \right)}

\newcommand{\supp}{\mathrm{supp}}

\newcommand{\fvg}{\mathcal{F}\left(V, G \right)}

\newcommand{\tfwg}{\tilde{\cf}(W, G)}

\newcommand{\fcvfg}{\mathcal{F}_c\left(V, \fg \right)}

\newcommand{\kp}{W^{k, p}}

\newcommand{\ka}{C^{k, \alpha}}

\newcommand{\fvr}{\mathcal{F}\left(V, \mathbb{R} \right)}

\newcommand{\fcvg}{\mathcal{F}_{c}\left(V, G \right)}

\newcommand{\fczvg}{\mathcal{F}^0_{c}\left(V, G \right)}

\newcommand{\dimv}{\dim_{\br} V}

\newcommand{\zero}{\mathbf{0}}
\newcommand{\one}{\mathbf{1}}

\newcommand{\vb}{\bar{v}}

\newcommand{\fttn}{\cf(\mathbb{T}^n, \mathbb{T}^n)}

\newcommand{\ta}{\tilde{a}}

\newcommand{\xv}{\mathfrak{X}_c(V)}
\newcommand{\xw}{\mathfrak{X}_c(W)}

\newcommand{\ty}{\tilde{y}}

\begin{document}

\title{A representation of range decreasing group homomorphisms}

\author{Ning Zhang}

\address{School of Mathematical Sciences \\ Ocean University of China \\Qingdao, 266100, P. R. China}

\email{nzhang@ouc.edu.cn}

\author{Lifan Liu}

\address{School of Mathematical Sciences \\ Ocean University of China \\Qingdao, 266100, P. R. China}

\email{21231111008@stu.ouc.edu.cn}



\begin{abstract}
The method of range decreasing group homomorphisms was first
introduced in \cite{z23} and later generalized in \cite{z25}.
This approach has proved effective for analyzing a broad class of weighted composition operators
between mapping spaces,
 including
  holomorphic maps, group homomorphisms, linear maps, semigroup
homomorphisms, Lie algebra homomorphisms
and algebra homomorphisms.
Previous work on range decreasing group homomorphisms has
focused mainly on specific subsets of mapping groups.
In this paper, we characterize range decreasing group homomorphisms on the entire mapping group.
As applications, we identify all such group homomorphisms
on certain mapping groups and classify
 a particular type of homomorphisms
between mapping groups.
 \end{abstract}

\maketitle

\pagestyle{myheadings} \markboth{\centerline{N. Zhang, L. Liu}}{\centerline{A representation of range decreasing group homomorphisms}}

\section{Introduction}

Throughout this paper, we assume that
all manifolds are
Hausdorff and all locally convex spaces are sequentially complete.
Unless otherwise indicated, we denote by $V$ a positive dimensional $C^{\infty}$ manifold,
possibly with boundary, that need not be paracompact or connected,
and by $G$ a positive or infinite dimensional locally exponential connected Lie group modelled
on a locally convex space (Notably, all Banach Lie groups fall within this category of locally exponential groups).
We fix a smoothness class $\cf =C^k_{MB}$
($C^k$ in the Michal-Bastiani sense, see \cite[Section I.2]{towards}),
$k=0, 1, \cdots, \infty$.
If $G$ is a Banach Lie group, we may also take $\cf$ to be
$C^k$ in the sense of Fr\'echet
differentiability for the same range of $k$, see \cite[Section 2.3]{amr}. If
$G$ is positive dimensional, we allow $\cf$ to be
locally H\"older $\ka$, $k=0, 1, \cdots$, $0 \le \alpha \le 1$,
where $C^{k, 0}=C^k$, or locally Sobolev $\kp$, $k=1, 2, \cdots$, $1 \le p <\infty$, $kp>\dimv$, or $k = \dimv$ and $p=1$.
The space $\fvg \subset C(V, G)$ of all $\cf$ maps $V \to G$
is a group under pointwise group operation.
If $V$ is compact, then $\fvg$ is a Lie group, see \cite[Theorem II.2.8]{towards} and \cite[Section 4(G)]{kr}.
Mapping spaces such as $\fvg$ and the operators between them
are of fundamental importance in analysis, geometry, mathematical physics and representation theory.

Suppose that $V$ and $W$ are finite dimensional $C^{\infty}$ manifolds, possibly
  with boundary, where $\dim_{\br} V \ge 1$ and
  $\dim_{\br} W \ge 0$,
  $\cf$ and $\tilde{\cf}$ are two smoothness classes,
  $H$ is a subgroup of $\fvg$ and
  $f\colon H \to \tfwg$ is a group homomorphism.
We say that $f$ is a weighted composition operator
if there exist maps $\phi\colon W \to V$
and $\gamma\colon W \times G \to G$ such that
 $\gamma(w, \cdot)\colon G \to G$ is a group isomorphism for each $w \in W$ and
 \begin{equation*}
 f(x)(w)=\gamma(w, x \circ \phi(w)), \hspace{2mm} w \in W, \hspace{2mm} x \in H.
 \end{equation*}
 If each $\gamma(w, \cdot)$ is the identity map on $G$, we call $f$ a composition or pullback operator.
 The constant maps form a subgroup of $\tfwg$ that
 can be identified with $G$.
 The composition operator induced by a constant map $\phi\colon W \to V$
 reduces to the evaluation map
 $$E_{\vb}\colon \fvg \ni x \mapsto x(\vb) \in G$$
 at a fixed point $\vb \in V$.

Weighted composition operators have been studied for many years, with various special cases investigated in the literature.
Let $\bff=\br$ or $\bc$. In \cite{g, mr}, Grabowski and Mr\v{c}un proved that
every algebra isomorphism $C^{\infty}(V, \bff) \to C^{\infty}(W, \bff)$
is the pullback by a $C^{\infty}$ diffeomorphism $W \to V$, extending earlier results in \cite{gk, myers, pursell}.
Let $\xv$ be the Lie algebra of all compactly supported $C^{\infty}$ vector fields on $V$.
The Shanks-Pursell theorem states that any Lie algebra isomorphism $f\colon$ $\xv$ $\to$ $\xw$
is of the form $f = d\psi$ for
 a $C^{\infty}$ diffeomorphism $\psi\colon V \to W$ \cite{sp, kmo}.
Consequently, $f$ can be regarded as a weighted composition operator in a certain sense, see \cite[(1.1)]{z25}.
For further results of Shanks-Pursell type, see \cite[Section X]{o} and \cite{gim, ab}.
Mr\v{c}un and \v{S}emrl established that
every isomorphism $C^k(V, \mathbb{R}) \to C^k(W, \mathbb{R})$ of multiplicative semigroups, for
$k \in \bn$,
is the pullback by a $C^k$ diffeomorphism \cite{ms}, generalizing
the earlier work by Milgram \cite[Theorem A]{mi}.
Sufficient conditions for a homomorphism $f$ between mapping groups
to be a weighted composition
operator were given
in \cite[Theorem 3.7]{hr} and \cite[Theorem 6.9]{fe}.
Notably, it follows from \cite[Theorem 4.1]{hr}
that for compact Hausdorff spaces $X$ and $Y$, any
linear bijection $f\colon C(X) \to C(Y)$ satisfying
$$f(C(X, \br \setminus \{0\}))=C(Y, \br \setminus \{0\})$$
must be a weighted composition operator.

Let $X$ be a subset of $\fvg$.
A map $f\colon X \to \tfwg$ is called range decreasing
if $f(x)(W) \subset x(V)$ for every $x \in X$.
Every composition operator $\fvg \ni x \mapsto x \circ \phi \in \tfwg$ is automatically range decreasing.
However,
a linear range decreasing map $\fvr \to \tilde{\cf}(W, \br)$ need not be
a composition operator \cite[P. 2190]{z23}.
Moreover,
for any $n \in \bn$,
one can construct a $C^{\infty}$ range decreasing map $C^{\infty}(S^2, \bc^n) \to C^{\infty}(S^2, \bc^n)$
that is not a composition operator \cite[P. 2180]{z23}.

In \cite[Theorem 1.1]{z23},
it was proved that
if $V$ is compact and connected and $2 \le \dim_{\br} G <\infty$,
then any range decreasing group homomorphism $f \colon \fvg \to G$
coincides with the evaluation map $E_{\vb}$ for some $\vb \in V$ on every
connected component of $\fvg$ that contains an element whose image is nowhere dense in $G$.
This result applies to the study of holomorphic maps between
mapping spaces. As another application, sufficient conditions were given in  \cite[Corollary 1.2]{z23}
for a group homomorphism $f\colon \cf^0(V, G) \to \tfwg$, where $\cf^0(V, G)$
is the identity component of $\fvg$,
to be a weighted composition operator.

According to \cite{z25}, analogues of the results in \cite[Theorem 1.1]{z23} hold
under much more general conditions.
We say that $x \in \fvg$ has compact support if $x = \one \in G$ outside a compact subset of $V$.
The subgroup of $\fvg$ consisting of all maps with compact support is denoted by $\fcvg$.
We define
 $\fczvg \subset \fvg$ as the subgroup of those $x$ for which
there exist a compact subset $K \subset V$ with $\supp x \subset K$
  and a homotopy $H\colon [0, 1] \times V \to G$
  relative to $V \setminus K$ such that $H(0, \cdot)$ is the constant map $\one$, $H(1, \cdot)=x$ and
  $H(t, \cdot) \in \fvg$ for all $t \in [0, 1]$.
  It is straightforward to verify that $\fczvg$ is a normal subgroup of $\fvg$.
  For $x \in \fvg$,  we denote by $[x]$ the coset of $\fczvg$ in $\fvg$
  containing $x$.
  If $2 \le \dim_{\br} G \le \infty$ and $f\colon \fvg \to G$ (respectively $f\colon \fcvg \to G$)
is a range decreasing group homomorphism with $f|_{\fczvg} \not\equiv \one$,
then there exists $\vb \in V$ such that $f(x)=x(\vb)$ for every $x$ in any
coset $[x_0]$ with
$x_0(V) \not=G$   \cite[Theorems 3.5 and 5.1]{z25}.
If $G$ is a locally convex space $\sce$ (in which case
$\cf_c^0(V, \sce)=\cf_c(V, \sce)$), then $f=E_{\vb}$ on its entire domain
\cite[Theorem 4.1 and Lemma 3.4]{z25}.
These results shed new light on several problems related to
 weighted composition operators. For instance, it has been
proved that the algebraic structure
of the space of smooth sections of an algebra bundle, whose typical
fiber is a positive dimensional simple unital algebra, completely determines
the bundle structure \cite[Theorem 12.1]{z25}.
When the algebra bundles are trivial with typical fiber $\bff=\br$ or $\bc$,
one recovers the results of Grabowski and Mr\v{c}un
 in \cite{g, mr}.
The Shanks-Pursell theorem has been extended from Lie algebra isomorphisms to Lie algebra homomorphisms
\cite[Theorem 11.1]{z25}.
The results of Mr\v{c}un and \v{S}emrl on isomorphisms
 of multiplicative semigroups
have been generalized to homomorphisms of multiplicative semigroups \cite[Theorem 9.5]{z25}.
A generalization of \cite[Theorem 4.1]{hr} appears in \cite[Theorem 8.1]{z25}.
Further applications include necessary and sufficient conditions for
a group isomorphism
between mapping groups to be a weighted composition operator \cite[Theorem 7.3]{z25}.

The method of range decreasing group homomorphisms is independent of the topology of the mapping spaces under consideration.
This approach applies to group homomorphisms in general, rather than being restricted to group isomorphisms or linear maps.
It can be used on any finite dimensional manifold, regardless of its paracompactness or connectedness,
and is compatible with a broad range of smoothness classes.

A general range decreasing group homomorphism
$f\colon \fvg \to \tfwg$ can be analyzed via
the family of homomorphisms $E_w \circ f\colon \fvg \to G$, $w \in W$.
Note that there exist range decreasing group homomorphisms
$C^{\infty}(S^3, \mathrm{SU}(2))$ $\to$ $\mathrm{SU}(2)$ that
do not take the form $E_{\vb}$
on certain components of $C^{\infty}(S^3, \mathrm{SU}(2))$
\cite[Proposition 4.2]{z25}.
In this paper, we characterize range decreasing group
homomorphisms $f\colon \fvg \to G$ (respectively $f\colon \fcvg \to G$)
 on their entire domain (Theorem \ref{general}).
 Together with the results of \cite{z25}, this yields
 a complete classification of range decreasing group homomorphisms
 $f\colon \fvg \to G$ with $f|_{\fczvg} \not\equiv \one$ (Corollary \ref{main}).
As applications, we
determine all range decreasing group homomorphisms
on certain mapping groups (Corollary \ref{compact} and
Theorem \ref{torus}). Furthermore,
 building on Theorem \ref{general}, we give a representation
 of a particular class of group homomorphisms between mapping groups (Corollary \ref{wco}),
 which generalizes \cite[Theorem 4.1]{hr} and \cite[Corollary 1.2]{z23}.
 In particular, by introducing new ideas into the proof of Corollary \ref{wco},
 we obtain a smoothness result similar to that of \cite[Corollary 1.2]{z23},
 but now for a broader class of
 manifolds $V$.

\section{Background and results}

A Lie group $G$ is called locally exponential
if it admits a $C^{\infty}$ exponential map $\exp _G\colon \fg \to G$, where $\fg$ is the Lie algebra of $G$,
that maps an open neighborhood
 of $\zero \in \fg$ diffeomorphically onto an open neighborhood of $\one \in G$
\cite[Section IV.1]{towards}.
For the topological degree of continuous maps between compact connected oriented manifolds
of the same dimension,
we refer to \cite[Section 5.1]{hi}.

\begin{thm} \label{general}
  Let  $V$ be a positive dimensional $C^{\infty}$ manifold,
possibly with boundary, that need not be paracompact or connected,
$G$ a positive or infinite dimensional locally exponential connected Lie group modelled
on a locally convex space,
  $f\colon \fvg \to G$ (respectively $f\colon \fcvg \to G$) a group homomorphism and
$Z(G)$ the center of $G$.
  If
  there exists $\vb \in V$ such that
  $f=E_{\vb}$ on the normal subgroup $\fczvg$, then
  there exists a group homomorphism $\psi\colon$ $\fvg/\fczvg$ $\to$ $Z(G)$
  (respectively $\psi\colon \fcvg/\fczvg \to Z(G)$) such that
    \begin{equation} \label{iff}
  f(x)= \psi([x]) x(\vb), \hspace{3mm} x \in \fvg \hspace{1mm} (\text{respectively} \hspace{2mm} x \in \fcvg),
  \end{equation}
  where $[x]$ denotes the coset of $\fczvg$
  containing $x$.
\end{thm}

From Theorem \ref{general} we immediately obtain the following

\begin{cor} \label{main}
  Let $V$ and $G$ be as in Theorem \ref{general}, and let $f\colon \fvg \to G$ (respectively $f\colon \fcvg \to G$) be a map.
  Then the following statements
  are equivalent.
  \begin{itemize}
    \item[(a)] The map $f$ is a group homomorphism, and there exists $\vb \in V$
    such that  $f(x)=x(\vb)$ for all $x$
    in every coset $[x_0]$ with $x_0$ non-surjective.
    In this case, $f$ is automatically range decreasing.

    \item[(b)] There exist a group homomorphism $\psi\colon$ $\fvg/\fczvg$ $\to$ $Z(G)$
    (respectively $\psi\colon \fcvg/\fczvg \to Z(G)$) and $\vb \in V$
  such that the kernel of $\psi$ contains
  the subset $\{[x_0]: x_0(V) \not=G \}$
  and  (\ref{iff}) holds.
  \end{itemize}
\end{cor}
 Assume that $\dim_{\mathbb{R}} G \geq 2$ and that Theorems 3.5 and 5.1 of \cite{z25} hold. Then
 any range decreasing group homomorphism $f\colon \fvg \to G$ (respectively $f\colon \fcvg \to G$)
  with $ f|_{\fczvg} \not\equiv \one$
   satisfies the conditions in Corollary \ref{main}(a).
   In this context, Corollary \ref{main}(b) provides a complete classification of such maps. Moreover,
 if $V$ is compact,
   the condition $ f|_{\fczvg} \not\equiv \one$
   is unnecessary for
   any range decreasing group homomorphism $f\colon \fvg \to G$, because
   every element of $G \subset \fczvg$ is a fixed point of $f$.

The following is an immediate consequence of \cite[Theorem 1.1]{z23} and Theorem \ref{general}.

\begin{cor} \label{compact}
  Let $V$ be a compact connected $C^{\infty}$ manifold, possibly with boundary,
  and let $G$ be a connected Lie group with $2 \le \dim_{\br} G <\infty$ and
  $Z(G)=\{\one\}$. Then every range decreasing group homomorphism $f\colon$ $\fvg$ $\to$ $G$
 is an evaluation map $E_{\vb}$ on its entire domain.
\end{cor}

Many connected Lie groups
have trivial center; examples of such groups include $\mathrm{SL}_{2k+1}(\br)$ and $\mathrm{SO}_{2k+1}(\br)$
for $k=1, 2, \cdots$ \cite[Example 9.3.13]{hn}.

Recall that if $V$ is compact, then $\fvg$ is a Lie group.
In this case, $\fczvg$ is the identity component of $\fvg$.
Furthermore, if $G$ is finite dimensional,
the
inclusion $\fvg$ $\hookrightarrow$ $C(V, G)$ is a homotopy equivalence, see
\cite[Theorem 13.14]{pa} and the remark following its proof.  Consequently,
the group of components of $\fvg$
is
$$\fvg/\fczvg \simeq \pi_0(\fvg) \simeq \pi_0(C(V, G)).$$
The Lie group $\mathrm{SU(2)}$ is diffeomorphic to the sphere $S^3$. We have
$Z(\mathrm{SU(2)}) \simeq \bz_2$ and $\pi_{10}(\mathrm{SU(2)}) \simeq \bz_{15}$ \cite[P. 76]{ntwo}.
There is a one-to-one correspondence between the free homotopy classes of continuous maps $S^{10} \to \mathrm{SU(2)}$
and the orbits of the action of $\pi_{1}(\mathrm{SU(2)}) \simeq \{\one\}$ on $\pi_{10}(\mathrm{SU(2)})$
\cite[Proposition 6.2.8]{t}.
Hence
$\pi_0(\cf(S^{10}, \mathrm{SU(2)}))$
consists of 15 elements.
This  implies that every group
homomorphism $\pi_0(\cf(S^{10}, \mathrm{SU(2)}))  \to Z(\mathrm{SU(2)})$
is constant $\one$,  even though $\pi_0(\cf(S^{10}, \mathrm{SU(2)}))$ $\not=$ $\{\one\}$ and $Z(\mathrm{SU(2)})$ $\not=$ $\{\one\}$.
By \cite[Theorem 1.1]{z23} and Theorem \ref{general}, every range decreasing
group homomorphism $\cf(S^{10}, \mathrm{SU(2)}) \to \mathrm{SU(2)}$ is of the form $E_{\vb}$
on its entire domain.

The homomorphism $\psi$ in Corollary \ref{main}(b) can be forced to be constant $\one$ even when
$G$ is commutative and $\fvg/\fczvg$ is infinite.

\begin{thm} \label{torus}
Let $\bt^n$ be the $n$ dimensional real torus, where $n \in \bn$.
The group $\pi_0(\fttn)$
is isomorphic to the additive group $M_n(\bz)$
  of $n \times n$ integer matrices, and the following statements hold.
\begin{itemize}
  \item[(a)] For $n \ge 2$, the group $\pi_0(\fttn)$ is generated by the
  subset $\{[x_0]: x_0(\bt^n) \not=\bt^n \}$, and
  every range decreasing group homomorphism $f\colon$ $\fttn$ $\to$ $\bt^n$
  is the evaluation map $E_{\vb}$
  at some $\vb \in \bt^n$ on its entire domain.

  \item[(b)] Let $f\colon$ $\cf(S^1, S^1)$ $\to$ $S^1$ be a range decreasing group homomorphism
  whose restriction to the subgroup $\cf_c^0(S^1, S^1)$ is an evaluation map.
  Then $f$ can be expressed as
  $$\cf(S^1, S^1) \ni x \mapsto z_0^{d(x)} x(\vb) \in S^1,$$
where $\vb, z_0 \in S^1$ are constants and $d(x)$ denotes the topological degree of $x$.
\end{itemize}
\end{thm}

For a Lie group $G$, we denote by $\autag$ the group of
algebraic group automorphisms of $G$ and by $\autg$ the group of Lie group automorphisms of $G$.
If $G$ is finite dimensional and connected,
then $\autg$ is a finite dimensional Lie group \cite[Theorem 2]{ho}.
As an application of Theorem \ref{general},  we generalize both \cite[Theorem 4.1]{hr}
and \cite[Corollary 1.2]{z23}.

\begin{cor} \label{wco}
  Let $V$ and $W$ be finite dimensional manifolds, possibly
  with boundary, with $\dim_{\br} V \ge 1$ and
  $\dim_{\br} W \ge 0$. Let
  $\cf$ and $\tilde{\cf}$ be smoothness classes, and let $f\colon \fvg \to \tilde{\cf}(W, G)$
  be a group homomorphism.
  Assume that every range decreasing group homomorphism $g\colon \fvg \to G$
  with $g|_{\fczvg} \not\equiv \one$
  satisfies the conditions in Corollary \ref{main}(a) (see \cite[Theorems 3.5 and 5.1]{z25}).
  Then the following statements are equivalent.
  \begin{itemize}
    \item[(a)] For every $w \in W$, the map $E_w \circ f|_{\fczvg} \not\equiv \one$
    and the map $E_w \circ f|_G\colon G \to G$
  is surjective. Furthermore, the following condition holds:
  \begin{equation}
  \label{nonzero}
    f(\cf(V, G \setminus \{\one\})) \subset \tilde{\cf}(W, G \setminus \{\one\}).
\end{equation}

    \item[(b)] There exist maps $\phi\colon W \to V$, $\gamma\colon W \to \autag$ and
    $$\psi_W\colon \fvg/\fczvg \times W \to Z(G)$$
    such that $\psi_W(\cdot, w)$ is a group
  homomorphism with $\{[x_0]: x_0(V) \not=G \} \subset \ker \psi_W(\cdot, w)$
  for every $w \in W$, and
  \begin{equation} \label{fxw}
    f(x)(w)=\gamma(w)(\psi_W([x], w) x \circ \phi(w)), \hspace{2mm} x \in \fvg, \hspace{2mm} w \in W.
  \end{equation}
  \end{itemize}
  Moreover, if $\dim_{\br} G <\infty$ and $E_w \circ f|_G \in \autg$ for every $w \in W$, then
  $\gamma\colon W \to \autg$ is an $\tilde{\cf}$ map.
\end{cor}

For the smoothness of the map $\phi$ in Corollary \ref{wco}(b), see \cite[Lemma 6.4]{z25}.
By (\ref{nonzero}), the homomorphisms $E_w \circ f|_G$, $w \in W$,
are injective. If $G$ is finite dimensional and the homomorphisms $E_w \circ f|_G$, $w \in W$,
are Borel measurable, then they are continuous \cite{kl}.
Consequently, they are injective Lie group homomorphisms and therefore belong to $\autg$.
If $G$ is a compact connected semisimple Lie group,
then every group
homomorphism $G \to G$ is automatically continuous by the van der Waerden theorem \cite{v}.
In this case, the homomorphisms $E_w \circ f|_G$, $w \in W$, are automatically elements of $\autg$.
For further details on the automatic continuity of group homomorphisms between Lie groups,
we refer to \cite{braun}.


\section{Proofs of the main results}

\noindent {\it Proof of Theorem \ref{general}. }
  Let $x_1 \in \fvg$ (respectively $x_1 \in \fcvg$).
  For every $x$ in the coset $[x_1]$, there exists a unique $y_{r}=y_r(x, x_1) \in \fczvg$
  such that $x=x_1 y_{r}$. Similarly, there exists a unique
  $y_{l} \in \fczvg$ such that $x=y_{l} x_1$.
  It follows that
  $$f(x)=f(x_1) y_{r}(\vb)=x_1(\vb) y_{r}(\vb) x_1(\vb)^{-1}f(x_1).$$
  Hence $x_1(\vb)^{-1} f(x_1)$ commutes with $y_{r}(\vb)$.

  Let $\exp _G\colon \fg \to G$ be the exponential map of $G$. For any $\tilde{a} \in \fg$,
  there exists
  $\tilde{y} \in \fcvfg$ such that $\tilde{y}(\vb)=\tilde{a}$. Note that
  $\exp_G \circ \tilde{y} \in \fczvg$. Since $\exp_G(\fg)$
  contains an open neighborhood of $\one \in G$, it generates the entire group $G$ \cite[Theorem 7.4]{hewitt}.
  Thus $x_1(\vb)^{-1} f(x_1) \in Z(G)$. Define a group homomorphism $\psi_1\colon$
  $\fvg$ $\to$ $Z(G)$ (respectively $\fcvg$ $\to$ $Z(G)$)
  by $\psi_1(x_1)=x_1(\vb)^{-1} f(x_1)$.
  Then $f=\psi_1 E_{\vb}$.
  Since $\psi_1|_{\fczvg} \equiv \one$,
  $\psi_1$ induces a homomorphism $\psi\colon$
  $\fvg/\fczvg$ $\to$ $Z(G)$
  (respectively $\psi\colon$ $\fcvg/\fczvg$ $\to$ $Z(G)$)
  such that (\ref{iff}) holds.
\qed

\vspace{2mm}

\noindent {\it Proof of Theorem \ref{torus}. }
Consider $\bt^n$ as the product $\prod_{i=1}^n S^1$. For $i=1, \cdots, n$,
let $\bt^{n, i} \subset \bt^n$ be the
subgroup of the form $\prod_{k=1}^n H_k$,
where $H_i=S^1$ and $H_k=\{1\} \subset S^1$ for $k \not=i$.
Let $P_j\colon \bt^n \to S^1$ be the projection onto
the $j$-th component of $\bt^n$, $j=1, \cdots, n$. For each $x \in \fttn$ and $i, j=1, \cdots, n$,
define maps $\xi_{ij, x}=P_j \circ x|_{\bt^{n, i}}\colon \bt^{n, i} \to S^1$.
  Then
  \begin{equation} \label{x}
  x(s_1, \cdots, s_n)=\te \left(\prod_{i=1}^n \xi_{i1, x}(s_i), \prod_{i=1}^n \xi_{i2, x}(s_i), \cdots, \prod_{i=1}^n \xi_{in, x}(s_i) \right),
  \end{equation}
  where $(s_1, \cdots, s_n) \in \bt^n$.
  Denote the topological degree
  of $\xi_{ij, x}$ by $d_{ij}(x)$.
  The matrix $D(x)=(d_{ij}(x)) \in M_n(\bz)$ depends only
  on the coset $[x] \in \pi_0(\fttn)$.

  Given $x_1, x_2 \in \fttn$ with  $D(x_1)=D(x_2)$,
  it follows from the Hopf degree theorem (\cite[Theorem 5.1.10]{hi}) that
  the maps $\xi_{ij, x_1}$ and $\xi_{ij, x_2}$
  are homotopic for all $i, j=1, \cdots, n$. In view of (\ref{x}), we have $[x_1]=[x_2]$.
  For any $A=(a_{ij}) \in M_n(\bz)$, define a Lie
  group homomorphism $x_A\colon \bt^n \to \bt^n$ by
  $$x_A(s_1, \cdots, s_n)=\te \left( \prod_{i=1}^n s_i^{a_{i1}}, \prod_{i=1}^n s_i^{a_{i2}}, \cdots, \prod_{i=1}^n s_i^{a_{in}} \right). $$
  Note that $D(x_A)=A$ and $x_A x_B=x_{A+B}$ for all $A, B \in M_n(\bz)$.
  Therefore
  $$\pi_0(\fttn) \simeq M_n(\bz).$$

  (a) For $i_0, j_0=1, \cdots, n$,
  let $E_{i_0 j_0}=(a_{ij}) \in M_n(\bz)$ be the matrix
   with $a_{i_0 j_0}=1$ and all other entries $0$.
  Note that the image of $x_{E_{i_0 j_0}}$
  is nowhere dense in $\bt^n$. The group $\pi_0(\fttn)$ is generated by
  the subset $\{[x_{E_{i j}}]: i, j=1, \cdots, n \}$.
  It follows from \cite[Theorem 1.1]{z23}
  and Theorem \ref{general} that every range decreasing group homomorphism $f\colon$ $\fttn$ $\to$ $\bt^n$
  is of the form $E_{\vb}$ on its entire domain.

 (b) By the definition of the topological degree, the coset $[x]$
 consists of surjective maps whenever $d(x) \not =0$. Every group homomorphism
 $\psi\colon \bz \to S^1$ is of the form $\psi(d)=z_0^d$ for some constant $z_0 \in S^1$.
 The conclusion of (b) then follows from Theorem \ref{general}.
 \qed

 \vspace{2mm}

\noindent {\it Proof of Corollary \ref{wco}. }
We first show that (b) implies (a).
Note that $\cf(V, G \setminus \{\one\})$
consists of non-surjective maps, and
$\psi_W([x], w)=\one$ for all $w \in W$ and all non-surjective $x \in \fvg$.
It is clear that (a) holds.

Next we show that (a) leads to (b).
Define maps
\begin{eqnarray*}
& \gamma\colon W \ni w \mapsto E_w \circ f|_G \in \autag \hspace{1mm} \text{and} & \\
& h_w=\gamma(w)^{-1} \circ (E_w \circ f)\colon \fvg \to G,  \hspace{2mm} w \in W. &
\end{eqnarray*}
Then $h_w|_{G}=\mathrm{id}$, $h_w|_{\fczvg} \not\equiv \one$
and $h_w(\cf(V, G \setminus \{\one\})) \subset G \setminus \{\one\}$.
For any $x \in \fvg$, we have $\di h_w\left(x (h_w(x))^{-1} \right)=\one$,
where $\di (h_w(x))^{-1} \in G \subset \fvg$ is the constant map taking this value.
Thus $\one \in x (h_w(x))^{-1}(V)$, which yields $h_w(x) \in x(V)$; that is, $h_w$ is range decreasing
for every $w \in W$.
Applying Corollary \ref{main}
to the homomorphisms $h_w$, $w \in W$, we obtain maps $\phi\colon W \to V$
  and
  $\psi_W\colon \fvg/\fczvg \times W \to Z(G)$
   such that $\psi_W(\cdot, w)$ is a group
  homomorphism with $\{[x_0]\colon x_0(V) \not=G \} \subset \ker \psi_W(\cdot, w)$ for every $w \in W$, and
  (\ref{fxw}) holds.

Finally, we consider the case where
 $\dim_{\br} G <\infty$ and $\gamma(w)=E_w \circ f|_G \in \autg$ for every $w \in W$.
Let $\tilde{\gamma}(w) \in Aut(\fg)$
be the Lie algebra automorphism induced by $\gamma(w)$, where $Aut(\fg)$  denotes the automorphism
group of the Lie algebra $\fg$ of $G$. Let $\exp_G\colon \fg \to G$
be the exponential map,
and choose an open ball $B$ centered at $\zero \in \fg$
such that $\exp_G|_{B}\colon B \to \exp_G(B)$ is a diffeomorphism.

For any $\ta_0 \in B \setminus \{\zero\}$ and any $w_0 \in W$,
we aim to prove that $\tilde{\gamma}(w)(\ta_0)$
is an $\tilde{\cf}$ map with respect to $w$ in a neighborhood of $w_0$.
We may assume that
$\tilde{\gamma}(w_0)(\ta_0) \in \frac{1}{2} B \setminus \{\zero\}$,
otherwise replace $\ta_0$ by a suitable scalar multiple.
Note that
\begin{eqnarray}
  & f \circ \exp_G(t\ta_0)(w_0) =  \gamma(w_0) \circ \exp_G(t\ta_0) & \notag \\
\label{comm} &  =\exp_G \circ \tilde{\gamma}(w_0)(t\ta_0) \in \exp_G(\te \frac{1}{2} B), \hspace{1mm} t \in [0, 1]. &
\end{eqnarray}
Choose a precompact
open connected neighborhood $O$ of $w_0$ such that
$$f \circ \exp_G(\ta_0)(w) \in \exp_G(\te \frac{1}{2} B) \setminus \{\one\}, \hspace{1mm} w \in \overline{O}.$$
Define a map
$\ty_0\colon \overline{O} \to \frac{1}{2} B \setminus \{\zero\}$
by
$$\ty_0(w)=\exp_G^{-1}(f \circ \exp_G(\ta_0)(w)), \hspace{1mm} w \in \overline{O}.$$
Then $\ty_0$ is an $\tilde{\cf}$ map on $O$.
For any $m \in \bn$, define two maps
\begin{eqnarray*}
& y_m =\exp_G \circ (\frac{\ty_0}{m})\colon \overline{O} \to \exp_G(\frac{1}{2m} B) \setminus \{\one\} \hspace{2mm} \text{and} &\\
& y^{\ast}_m =f \circ \exp_G(\frac{\ta_0}{m}) \in \tfwg. &
\end{eqnarray*}
Let $O_m \subset O$ be the closed subset $\{w \in O: y_m(w)=y^{\ast}_m(w)\}$.
Then $O_1=O$.
By (\ref{comm}), we have $\ty_0(w_0)=\tilde{\gamma}(w_0)(\ta_0)$, which implies that
$w_0 \in O_m$ for every $m \in \bn$.
Let $w_1 \in O_m$. Then $\ty_0(w_1)+m \ta \in B$ for all $\ta \in \frac{1}{2m} B$. Consequently,
the map
\begin{eqnarray*}
  & \exp_G\left(\frac{\ty_0(w_1)}{m}+\frac{1}{2m} B\right) \ni \exp_G\left(\frac{\ty_0(w_1)}{m}+\ta \right)&\\
  & \mapsto \left(\exp_G\left(\frac{\ty_0(w_1)}{m}+\ta \right)\right)^m=\exp_G(\ty_0(w_1)+m \ta) \in \exp_G(B) &
\end{eqnarray*}
 is injective.
Hence, for every $w \in O$, the set $\exp_G\left(\frac{\ty_0(w_1)}{m}+\frac{1}{2m} B\right)$
contains at most one $m$-th root of $f \circ \exp_G(\ta_0)(w)$.
Note that
$$y_m^m(w)=f \circ \exp_G(\ta_0)(w)=(y^{\ast}_m)^m(w), \hspace{1mm} w \in O.$$
By continuity of $y_m$ and $y^{\ast}_m$,
$w_1$ is an interior point of $O_m$. Thus $O_m=O$ for every $m \in \bn$.
Since $y^{\ast}_m(w)=\gamma(w) \circ \exp_G(\frac{\ta_0}{m})$ for $w \in O$, we have
\begin{equation} \label{m}
  \exp_G\left(\te \frac{\ty_0(w)}{m}\right)=\exp_G\left(\te \frac{\tilde{\gamma}(w)(\ta_0)}{m}\right), w \in O, m \in \bn.
\end{equation}
For any $w \in O$,
the restriction of $\exp_G$ to the one dimensional
subspace $\{\ty_0(w) \br \} \subset \fg$ is a Lie
group homomorphism, and its kernel is either $\{0\}$ or a cyclic subgroup.
From (\ref{m}) we therefore conclude that
\begin{equation*}
\tilde{\gamma}(w)(\ta_0)=\ty_0(w), w \in O.
\end{equation*}

Thus the map $W \ni w \mapsto \tilde{\gamma}(w) \in Aut(\fg)$ can be interpreted
as a matrix valued $\tilde{\cf}$ map.
Recall that the map $\autg \ni \sigma \mapsto d_{\one} \sigma \in Aut(\fg)$, where $d_{\one} \sigma$
denotes the differential of $\sigma$ at the identity $\one$,
is an injective Lie group homomorphism onto a closed subgroup of $Aut(\fg)$ (see, e.g., \cite[Subsection 11.3.1]{hn}).
Consequently, $\gamma$ is an
$\tilde{\cf}$ map.
 \qed

\noindent * Corresponding author: Ning Zhang, nzhang@ouc.edu.cn, 
School of Mathematical Sciences, Ocean University of China, Qingdao, 266100, P. R. China.

\section*{Funding}

This research was partially supported by the Natural Science Foundation of Shandong Province of China
grant ZR2023MA035.

\section*{Declarations}

Consent to Participate declaration: not applicable.

Consent to Publish declaration: not applicable.

Ethics declaration: not applicable.

\end{document}